\documentclass{article}
\usepackage{latexsym, amsfonts, amsmath}
\newcommand{\be}{\begin{equation}}
\newcommand{\ee}{\end{equation}}
\newcommand{\bea}{\begin{eqnarray}}
\newcommand{\eea}{\end{eqnarray}}
\newcommand{\bean}{\begin{eqnarray*}}
\newcommand{\eean}{\end{eqnarray*}}

\newcommand{\brray}{\begin{array}}
\newcommand{\erray}{\end{array}}
\newcommand{\ben}{\begin{equation}{nonumber}}
\newcommand{\een}{\end{equation}{nonumber}}

\newtheorem{dfn}{Definition}[section]
\newtheorem{thm}[dfn]{Theorem}
\newtheorem{lmma}[dfn]{Lemma}
\newtheorem{ppsn}[dfn]{Proposition}
\newtheorem{crlre}[dfn]{Corollary}
\newtheorem{xmpl}[dfn]{Example}
\newtheorem{rmrk}[dfn]{Remark}
\newcommand{\bdfn}{\begin{dfn}}
\newcommand{\bthm}{\begin{thm}}
\newcommand{\blmma}{\begin{lmma}}
\newcommand{\bppsn}{\begin{ppsn}}
\newcommand{\bcrlre}{\begin{crlre}}
\newcommand{\bxmpl}{\begin{xmpl}}
\newcommand{\brmrk}{\begin{rmrk}}
\newcommand{\edfn}{\end{dfn}}
\newcommand{\ethm}{\end{thm}}
\newcommand{\elmma}{\end{lmma}}
\newcommand{\eppsn}{\end{ppsn}}
\newcommand{\ecrlre}{\end{crlre}}
\newcommand{\exmpl}{\end{xmpl}}
\newcommand{\ermrk}{\end{rmrk}}

\newcommand{\IC}{\mathbb{C}}

\newcommand{\cla}{{\cal A}}

\newcommand{\clc}{{\cal C}}

\newcommand{\cll}{{\cal L}}

\newcommand{\clq}{{\cal Q}}

\def\a*{{\cal A}_{h,*}}
\def\B{{\cal B}(h)}
\def\B1{{\cal B}_1(h)}
\def\b{{\cal B}^{\rm s.a.}(h)}
\def\b1{{\cal B}^{\rm s.a.}_1(h)}

\newcommand{\ot}{\otimes}

\newcommand{\raro}{\rightarrow}

\def \qed {$\Box$}

\def\a*{{\cal A}_{h,*}}
\def\B{{\cal B}(h)}
\def\B1{{\cal B}_1(h)}
\def\b{{\cal B}^{\rm s.a.}(h)}
\def\b1{{\cal B}^{\rm s.a.}_1(h)}

\begin{document}

\begin{center}
{\Large{\bf An averaging trick for smooth actions of compact quantum groups on manifolds}}\\ 
{\large  \bf Debashish Goswami \footnote{Partially supported by Swarnajayanti Fellowship from D.S.T. (Govt. of India)}, 
\bf Soumalya Joardar \footnote{Acknowledges support from CSIR}}\\ 
Indian Statistical Institute\\
203, B. T. Road, Kolkata 700108\\
Email: goswamid@isical.ac.in\\
\end{center}
\begin{abstract}
We prove that, given  any smooth action of a compact quantum group (in the sense of \cite{rigidity}) on a compact smooth manifold satisfying 
some more natural conditions, one can get a Riemannian structure on the manifold for which the corresponding $C^\infty(M)$-valued inner product on the space of one-forms is 
 preserved by the action.

\end{abstract}
\begin{center}
 {\it Dedicated to Prof. Kalyan B. Sinha on his seventieth birthday}.
\end{center}

\section{Introduction}
It is both interesting as well as important to study quantum group actions on classical (commutative) and noncommutative spaces. 
Indeed, quantum group actions can be viewed as generalised symmetries of a classical or quantum system modelled by commutative or noncommutative manifolds.
In this context, it is natural ask the question whether one can have genuine (i.e. which are not groups) compact quantum group actions on (compact) 
classical spaces. Indeed, this has an affirmative answer in general. First examples of this kind were produced by S. Wang (\cite{wang}, 
see also later  of other mathematicians in this direction, e.g. \cite{ban}, \cite{ban2} etc.) who defined and studied a quantum-group generalisation of the group of permutations of $n$ objects, called the quantum permutation group, and gave its action on the algebra of functions on finite set of cardinality $n$.  For $n \geq 4$ this quantum group is a genuine one. However, in all such cases the underlying set is disconnected. It took quite a long 
time since the work of Wang before H. Huang (\cite{huichi}) came up with 
 several example of genuine compact quantum groups acting faithfully on compact connected topological spaces. On the other hand, there were indications (e.g. \cite{banica_jyotish}) 
  that such a construction 
 would not be possible if the space is a connected smooth manifold. One of the author of the present paper (D. Goswami) made this conjecture and both the authors could prove 
  (\cite{rigidity}) the non-existence of any faithful action of a genuine compact quantum group on a compact connected manifold 
 if the action is assumed to be smooth and isometric  in a natural sense. In this context, it turns out to be useful to prove an analogue of the classical averaging 
  technique for compact group actions on Riemannian manifolds. The aim of this note is to achieve such a result for a smooth action of a compact quantum group on 
   a compact Riemannian manifold under certain natural conditions which are valid for a large class of examples. We hope it has the potential of 
  generalisation to the context of noncommutative manifold a la Connes 
(\cite{Connes}).

\section{Notaton and preliminaries}
We denote by $\hat{\ot}$ spatial (minimal) $C^*$ tensor product of $C^*$ algebras.
\bdfn A compact quantum group (CQG for short) is a  unital $C^{\ast}$ algebra $\clq$ with a
coassociative coproduct 
 $\Delta$ from $\clq$ to $\clq \hat{\ot} \clq$  
such that each of the linear spans of $\Delta(\clq)(\clq\ot 1)$ and 
 $\Delta(\clq)(1\ot \clq)$ is norm-dense in $\clq \hat{\ot} \clq$. 
 
An action of $\clq$ on a unital $C^*$ algebra $\clc$ is a unital $\ast$-homomorphism  $\alpha:\clc\raro \clc \hat{\ot}\clq$ such that 
 $(\alpha \ot id)\alpha=(id\ot \Delta)\alpha$ and 
  $\overline{Sp} \{  \alpha(\clc)(1\ot \clq)\} =\clc\hat{\ot}\clq$.
  \edfn
   We denote by $\clq_0$ the dense unital Hopf $\ast$ algebra in $\clq$ spanned by the matrix coeeficients of irreducible unitary representations of $\clq$ (see, e.g. \cite{Van}).
  Given an action $\alpha$ of $\clq$ on $\clc$ we always get a dense unital $\ast$-subalgebra $\clc_0$ of $\clc$ on which $\alpha$ is algebraic, i.e. maps $\clc_0$ to the algebraic tensor 
   product $\clc_0 \ot \clq_0$.
  
\bdfn
An action $\alpha$ is said to be faithful if the $\ast$-subalgebra
of $\clq$ generated by \{ $(\omega\ot id)\alpha(a),~\omega \in \clc^*,~a \in \clc\}$, where
$\clc^*$ is the set of  bounded linear functionals on $\clc$, is dense in $\clq$.
\edfn

We refer \cite{Woronowicz}, \cite{Van} for the theory  of unitary representation of CQG's and to \cite{Connes} for the 
 framework of noncommutative geometry given by spectral triples. 
 
 \bdfn
 For a compact Riemannian manifold $M$, we say that action $\alpha$ of a CQG $\clq$ on $C(M)$ to be isometric  if it maps $C^\infty(M)$ to $C^\infty(M, \clq)$
  and for every bounded linear functional $\phi$ on $\clq$, one has $\cll \circ \alpha_\phi=\alpha_\phi \circ \cll$, where $\alpha_\phi=
   ({\rm id} \ot \phi)\circ \alpha$ and $\cll$ is the restriction of the Hodge Laplacian $-d^*d$ on the space of smooth functions. 
 \edfn
 The following result is proved in \cite{Goswami}.
 \bppsn 
 For a compact Riemannian manifold $M$, there is a universal object in the category of CQG's having isometric actions on $M$. We call this CQG the quantum
  isometry group of $M$.
 \eppsn
 
    
      An action $\alpha$ of a CQG $\clq$ on $C(M)$  is called  smooth if  it  maps 
 $C^\infty(M)$ to $C^\infty(M, \clq)$ and the span of $\alpha(C^\infty(M))(1 \ot \clq)$ is dense in $C^\infty(M, \clq)$ in the natural Frechet
  topology. It has been proved in \cite{rigidity} that smooth actions are automatically continuous as a map from $C^\infty(M)$ to $C^\infty(M, \clq)$ in the respective 
   Frechet topologies coming from that of $C^\infty(M)$. For any $C^*$-algebra $\clc$, We consider the set  of smooth $\clc$-valued one-forms $\Omega^1(M, \clc)$ with the natural 
  Frechet topology coming from $M$ ($\Omega^1(M) :=\Omega^1(M,\IC)$) and  the obvious  $C^\infty(M,\clc)$-bimodule structure. Given  a smooth action 
   $\alpha$ of $\clq$ we call a continuous $\IC$-linear map $\Gamma: \Omega^1(M) \raro \Omega^1(M, \clq)$ to be a representation if $\Gamma$ is 
    co-associative in the obvious sense and $\Gamma(\xi f)=\Gamma(\xi) \alpha(f)=\alpha(f) \Gamma(\xi)$ for $\xi \in \Omega^1(M), f \in C^\infty(M)$. 
  
  We often say that $\alpha$ is a smooth action on $M$ to mean that it is a smooth action on $C(M)$ in the sense discussed above.  
  For such an action we denote $(d \ot {\rm id}) (df)$ by $d\alpha(df)$. The $C^\infty(M)$-valued inner product on $\Omega^1(M)$ coming from the Riemannian structure is denoted by $<< \cdot, \cdot>>$ and we say that a smooth action $\alpha$ preserves the Riemannian structure (or the Riemannian inner product) if $<<d\alpha(df), d\alpha (dg)>>=\alpha(<<df,dg>>)$ for all $f,g \in C^\infty(M)$. It is proved in \cite{rigidity} that a smooth action on a compact Riemannian manifold $M$ (without boundary) is isometric if and only if it preserves the inner product.

  Before we state and prove the main result in the next section, let us collect a few facts about a smooth faithful action of compact quantum groups on compact manifolds, for the details 
   of which the reader may be referred to \cite{rigidity}  and references therein.
   \bppsn
   If a CQG $\clq$ acts faithfully and smoothly on a smooth compact manifold $M$ then we have:\\
   (i) $\clq$ has a tracial Haar state, i.e. it is Kac type CQG.\\
   (ii) The action is injective.\\
   (iii) The antipode $\kappa$ satisfies $\kappa(a^*)=\kappa(a)^*$.
   \eppsn
   
   We usually denote by $\otimes$ algebraic tensor product of vector spaces or algebras. We also use Sweedler convention for Hopf algebra coproduct as well as its analogue for (co)-actions 
    of Hopf algebras. That is, we simply write $\Delta(q)=q_{(1)} \ot q_{(2)}$ suppressing finite summation, where $\Delta$ denote the co-product map of a Hopf algebra and $q$ is an element 
   of the Hopf algebra. Similarly, for an algebraic (co)action $\alpha$ of a Hopf algebra on some algebra $\clc$, we write $\alpha(a)=a_{(0)} \ot a_{(1)}.$
\section{The main result}

Fix  a compact Riemannian 
manifold $M$ (not necessarily orientable) and a smooth action $\alpha$ of a CQG 
$\clq$. We make the following assumptions for the rest of thye paper.\\
{\bf Assumption I}: There is a Fr\'echet dense unital $\ast$-subalgebra $\cla$ of $C^\infty(M)$ such that 
$ << d\alpha(df), d \alpha(dg) >> \in \cla$ for all $f,g \in \cla$.\\
{\bf Assumption II}: There is a well-defined representation $\Gamma $ on $\Omega^1(M)$ in the sense discussed earlier, such that $\Gamma(df)=(d \ot {\rm id})(\alpha(f))$ for all $f \in C^\infty(M)$.
 We'll denote this $\Gamma$ by $d\alpha$.\\

We now state and prove the main result that we can equip $M$ with a new Riemannian
structure with respect to which the action becomes inner product preserving
using an analogue of the averaging
technique of classical differential geometry. 
\bthm
\label{average}
$M$ has a Riemannian structure such that $\alpha$ is inner product
preserving.\\
\ethm
 
 Note that the first assumption holds for a large class of examples, such as algebraic actions of CQG's compact, smooth, real varieties where 
  the complexifed coordinate algebra of the variety can be chosen as $\cla$. On the other hand, the second assumption means that the action on $M$ in some sense lifts to the space of one-forms. This is always automatic for a smooth  action by (not necessarily compact) 
  groups, and in fact is nothing 
  but the differential of the map giving the action. Moreover, it is easy to see that any CQG action which preserves the Riemannian inner product does admit 
   such a lift on the bimodule of one-forms, i.e. satisfies the assumption II. Therefore, it is a reasonable assumption too.

{\it Proof of Theorem \ref{average}}:\\
 
\indent We break the proof of into a number of lemmas.
\blmma
Define the following map $\Psi$ from $\cla \ot \clq_0$ to $\cla$ :
$$ \Psi(F):=({\rm id}\ot h)({\rm id } \ot m)({\rm id} \ot \kappa\ot {\rm id})(\alpha\ot
{\rm id})(F).   $$ Here $m:\clq_{0}\ot\clq_{0}\raro \clq_{0}$ is the multiplication map.
Then $\Psi$ is a completely positive map.
\elmma
{\it Proof:}\\
As the range is a subalgebra of a unital commutative $C^*$ algebra, it is enough to prove positivity. Let $F=G^*G$ in $\cla \ot \clq_0$ where $G=\sum_i f_i \ot q_i$, (finite sum)
 for some $f_i \in \cla, q_i \in \clq_0$. We write $\alpha(f)=f_{(0)} \ot f_{(1)}$ in Sweedler notation as usual, and observe that 
 \bean \lefteqn{\Psi(F)}\\
 & =& \sum_{ij}f_{i(0)}^* f_{j(0)}h(\kappa(f_{i(1)}^*f_{j(1)})q_i^*q_j)\\
 &=& \sum_{ij}f_{j(0)}f_{i(0)}^*h(q_j(\kappa(f_{j(1)}))^*\kappa(f_{i(1)})q_i^*)\\
 &=& ({\rm id} \ot h)(\xi^*\xi) \geq 0,\eean\\
 where $\xi=\sum_i f_{i(0)}^* \ot  \kappa(f_{i(1)})q_i^*,$ and note also that we have used above the facts that $h$ is tracial and $\kappa$ is $\ast$-preserving. \qed\\
 
For
$\omega,\eta\in
\Omega^{1}(\cla)$ We define 
$$<<\omega,\eta>>^{'}:=\Psi(<<d\alpha(\omega),d\alpha(\eta)>>),$$
 which is well defined as we have assumed that
$<<d\alpha(ds_{1}),d\alpha(ds_{2}))>>\in \cla\ot\clq_{0}$ for
$s_{1},s_{2}\in \cla$. Moreover, by complete positivity of $\Psi$ this gives a non-negative definite sesquilinear form on $\Omega^1(\cla)$. 
As the action is algebraic over $\cla$, we shall use Sweedler's notation to
prove the following 
\blmma
\label{new_inner}
For $\omega,\eta\in \Omega^{1}(\cla), \ f\in \cla$,
$<<\omega,\eta>>^{'}=(<<\eta,\omega>>^{'})^{ \ast}$ and $<<\omega,\eta
f>>^{'}=<<\omega,\eta>>^{'}f$
\elmma
{\it Proof}:\\
It is enough to prove the lemma for $\omega=d\phi$ and $\eta=d\psi$ for
$\phi,\psi \in \cla$. First observe that as we have $\kappa=\kappa^{-1}$, for
$z\in \clq_{0}$ applying $\kappa$ on $z_{(1)}\kappa(z_{(2)})=\epsilon(z).1$, we
get
\begin{eqnarray} 
z_{(2)}\kappa(z_{(1)})=\epsilon(z).1.
\end{eqnarray}
 We denote $<<d\phi_{(0)},d\psi_{(0)}>>$ by $x$ and
$\phi_{(1)}^{\ast}\psi_{(1)}$ by $y$.
Then \begin{eqnarray*}
&& <<d\phi,d\psi f>>^{'}\\
&=& (id \ot h)(id \ot m)(id \ot \kappa \ot id)(\alpha \ot
id)<<d\alpha(d\phi),d\alpha(d\psi f)>>\\
&=& (id \ot h)(id \ot m)(id \ot \kappa \ot id)(\alpha \ot id)(xf_{(0)}\ot
yf_{(1)})\\
&=& (id \ot h)(id \ot m)(id \ot \kappa \ot id) (x_{(0)}f_{(0)(0)} \ot
x_{(1)}f_{(0)(1)} \ot yf_{(1)})\\
&=& (id \ot h)(x_{(0)}f_{(0)(0)} \ot \kappa(x_{(1)}f_{(0)(1)})yf_{(1)})\\
&=& x_{(0)}f_{(0)(0)} h(f_{(1)}\kappa(f_{(0)(1)})\kappa(x_{(1)})y) ( by \
tracial \ property \ of \ h)\\
&=& x_{(0)}f_{(0)}h(f_{(1)(2)}\kappa (f_{(1)(1)})\kappa(x_{(1)})y)\\
&=& x_{(0)}f_{(0)}h(\epsilon(f_{(1)}).1.\kappa(x_{(1)})y) \\
&=& x_{(0)}(id\ot \epsilon)\alpha(f)h(\kappa(x_{(1)})y)\\
&=& x_{(0)}fh(\kappa(x_{(1)})y).
\end{eqnarray*}
On the other hand, \begin{eqnarray*}
<<d\phi,d\psi>>^{'}f &=& [(id\ot h)(id\ot m)(id\ot \kappa\ot id)(\alpha \ot
id)<<d\alpha(d\phi),d\alpha(d\psi)>>]f\\
&=& [(id\ot h) (id\ot m)(id\ot \kappa \ot id)(x_{(0)}\ot x_{(1)}\ot y)]f\\
&=& x_{(0)}fh(\kappa (x_{(1)})y).
\end{eqnarray*}
Also we have
\begin{eqnarray*}
 && <<d\phi,d\psi>>^{'}\\
 &=& (id \ot h)(id \ot m)(id \ot \kappa \ot id)(\alpha \ot
id)(<<d\phi_{(0)},d\psi_{(0)}>>\ot \phi_{(1)}^{\ast}\psi_{(1)})\\
 &=& (id \ot h)(id \ot m)(id \ot \kappa \ot id)(\alpha \ot
id)(<<d\psi_{(0)},d\phi_{(0)}>>^{\ast}\ot \phi_{(1)}^{\ast}\psi_{(1)})\\
 &=& (id \ot h)(id \ot m)(id \ot \kappa \ot
id)(<<d\psi_{(0)},d\phi_{(0)}>>_{(0)}^{\ast}\ot
<<d\psi_{(0)},d\phi_{(0)}>>_{(1)}^{\ast}\ot \phi_{(1)}^{\ast}\psi_{(1)})\\
 &=&<<d\psi_{(0)},d\phi_{(0)}>>^{\ast}h((\kappa(<<d\psi_{(0)},d\phi_{(0)}>>))^{
\ast} \phi_{(1)}^{\ast}\psi_{(1)}) ( \ since \ \kappa \ is \ \ast \
preserving)\\
\end{eqnarray*}
\indent Hence we have 

$$<<d\phi,d\psi>>^{' \ast}
= <<d\psi_{(0)},d\phi_{(0)}>>
h((\kappa(<<d\psi_{ (0)},d\phi_{(0)}>>))\psi_{(1)}^{\ast}\phi_{(1)})$$
$( \ since
\ h \ is \ tracial \ and \ h(a^{\ast})=\overline{h(a)})$.\\

But we can readily see that 
$$<<d\psi,d\phi>>^{'}=<<d\psi_{(0)},d\phi_{(0)}>>
h((\kappa(<<d\psi_{ (0)},d\phi_{(0)}>>))\psi_{(1)}^{\ast}\phi_{(1)}),$$
which completes the proof of the lemma.\qed\\
Actually we can extend $<<,>>^{'}$ to a slightly bigger set than
$\Omega^{1}(\cla)$ namely $\Omega^{1}(\cla)C^{\infty}(M)= \ Sp \ \{\omega
f: \omega\in \Omega^{1}(\cla),f\in C^{\infty}(M)\}$.\\
\indent For $\omega,\eta\in \Omega^{1}(\cla)C^{\infty}(M)$,
$\omega=\sum\omega_{i}f_{i} \ and \ \eta=\sum\eta_{i}g_{i}( \ finite \  sums)$,
$\omega_{i},\eta_{i}\in \Omega^{1}(\cla) \ and \ f_{i},g_{i}\in C^{\infty}(M)$
(say) we can choose sequences $f_{i}^{(n)},g_{i}^{(n)}$ from $\cla$ such that
$f_{i}^{(n)}\raro f_{i}$ and $g_{i}^{(n)}\raro g_{i}$ in the
corresponding Fr\'echet topology and by
Lemma \ref{new_inner} observe that 
\begin{eqnarray}
 &&<<\sum_{i} \omega_{i}f_{i}^{(n)},\sum_{j}\eta_{j}g_{j}^{(n)}>>^{'}
\nonumber \\
 &=& \sum_{i,j}
\overline{f_{i}^{(n)}}<<\omega_{i},\eta_{j}>>^{'}g_{j}^{(n)}\nonumber \\
 &\raro& \sum_{i,j}\overline{f_{i}}<<\omega_{i},\eta_{j}>>^{'}g_{j}:=
<<\omega,\eta>>^{'}
\end{eqnarray}
Clearly this definition is independent of the choice of sequences $f_{i}^{(n)}$
and $g_{i}^{(n)}$.\\
We next prove the following
\blmma
\label{new_inner_preserve}
For $\phi,\psi\in \cla$,
\begin{eqnarray}
<<d\alpha(d\phi),d\alpha(d\psi)>>^{'}= \alpha(<<d\phi,d\psi>>^{'})
\end{eqnarray}
\elmma
{\it Proof}:\\
With $x,y$ as before we have

 {\bf Claim 2}: We can extend the definition of $<<,>>^{'}$ for $\omega,\eta \in
\Omega^{1}(\cla)C^{\infty}(M)$ such that
 \begin{eqnarray}
 \label{1111}
  \forall \ \ f\in C^{\infty}(M), <<(d\phi),(d\psi)f>>^{'}=<<d\phi,
d\psi>>^{'}f 
  \end{eqnarray}
 {\it Proof}:\\
 For $f\in C^{\infty}(M)$, define $<<(d\phi),(d\psi)f>>^{'}:= \ 
lim<<d\phi, d\psi \ f_{n}>>^{'}$, where $f_{n}\in \cla$ with $lim \
f_{n}= f$, where the limits are taken in the Fr\'echet topology.\\
 Observe that $<<d\phi, d\psi \ f_{n}>>^{'}$ is Fr\'echet Cauchy as 
 \begin{eqnarray*}
 &&  <<d\phi, d\psi \ f_{n}>>^{'}- <<d\phi, d\psi \ f_{m}>>^{'}\\
 &=& <<d\phi,d\psi>>^{'}(f_{n}- f_{m})
 \end{eqnarray*}
 So $<<d\phi, d\psi  f>>^{'}= lim<<d\phi, d\psi >>^{'}f_{n}= <<d\phi,
d\psi
 >>^{'}f$, again the limit is taken in the corresponding Fr\'echet topology.\\
 That proves the claim.\\
\begin{eqnarray*}
&& <<d\alpha(d\phi),d\alpha (d\psi)>>^{'}\\
&=& (id\ot h\ot id)(id\ot m\ot id)(id\ot \kappa\ot id\ot id)(\alpha\ot id\ot
id)(x\ot \Delta(y))\\
&=& (id\ot h\ot id)(id\ot m\ot id)(id\ot \kappa\ot id\ot id)(x_{(0)}\ot
x_{(1)}\ot y_{(1)}\ot y_{(2)})\\
&=& (id\ot h\ot id)(x_{(1)}\ot \kappa(x_{(2)})y_{(1)}\ot y_{(2)})\\
&=& x_{(0)}\ot h(\kappa(x_{(1)})y_{(1)})y_{(2)}.
\end{eqnarray*}
\vspace{0.1in}
On the other hand 
\begin{eqnarray*}
\alpha(<<d\phi,d\psi>>^{'})&=& x_{(0)(0)}h(\kappa(x_{(1)})y)\ot x_{(0)(1)}\\
&=& x_{(0)}\ot x_{(1)(1)}h(\kappa(x_{(1)(2)})y)\\
&=& x_{(0)}\ot x_{(1)(1)}h(\kappa(y)(x_{(1)(2)})) ( \ since \
h(\kappa(a))=h(a))\\
\end{eqnarray*}  
Hence it is enough to show that $h(\kappa(c)b_{(2)})
b_{(1)}=h(\kappa(b)c_{(1)})c_{(2)}$ where $b,c\in \clq_{0}$, for then taking
$x_{(1)}=b$ and $y=c$ we can complete the proof.\\
\indent We make the transformation $T(a\ot b)=\Delta(\kappa(a))(1\ot b)$. \\
Then
\begin{eqnarray*}
&&(h\ot id)T(a\ot b)\\
&=& (h\ot id)\Delta(\kappa(a))(1\ot b)\\
&=& ((h\ot id)\Delta(\kappa(a)))b\\
&=& h(\kappa(a))b\\
&=& (h\ot id)(a\ot b)\\
\end{eqnarray*}
Hence $h(b_{(2)}\kappa(c))b_{(1)}=(h\ot id)T(b_{(2)}\kappa(c)\ot b_{(1)})$.\\
So, by using traciality of $h$ it is enough to show that $T(b_{(2)}\kappa(c)\ot
b_{(1)})= c_{(1)}\kappa(b)\ot c_{(2)}$.\\
\begin{eqnarray*}
&& T(b_{(2)}\kappa(c)\ot b_{(1)})\\
&=& \Delta(\kappa(b_{(2)}\kappa(c)))(1\ot b_{(1)})\\
&=& \Delta(c \kappa(b_{(2)}))(1\ot b_{(1)})\\
&=& (c_{(1)}\ot c_{(2)})[\kappa(b_{(2)(2)})\ot \kappa(b_{(2)(1)})](1\ot
b_{(1)})\\
&=& (c_{(1)}\ot c_{(2)})m_{23}(\kappa(b_{(2)(2)})\ot \kappa(b_{(2)(1)})\ot
b_{(1)})\\
&=& (c_{(1)}\ot c_{(2)})m_{23}(\kappa\ot \kappa\ot id)\sigma_{13}( b_{(1)}\ot
b_{(2)(1)}\ot b_{(2)(2)})\\
&=& (c_{(1)}\ot c_{(2)})m_{23}(\kappa\ot \kappa\ot id)\sigma_{13}(b_{(1)(1)}\ot
b_{(1)(2)}\ot b_{(2)})\\
&=& (c_{(1)}\ot c_{(2)})m_{23}(\kappa(b_{(2)}\ot \kappa (b_{(1)(2)})\ot
b_{(1)(1)})\\
&=& (c_{(1)}\ot c_{(2)})(\kappa(b_{(2)})\ot \epsilon(b_{(1)}).1_{\clq}) (by \
(10))\\
&=& (c_{(1)}\ot c_{(2)})(\kappa\ot \kappa)((b_{(2)})\ot
\epsilon(b_{(1)}).1_{\clq})\\
&=& (c_{(1)}\ot c_{(2)})(\kappa\ot \kappa)(\epsilon(b_{(1)})b_{(2)}\ot
1_{\clq})\\
&=&  c_{(1)}\kappa(b)\ot c_{(2)}
\end{eqnarray*}
Which proves the claim.\\
\vspace{0.2in} \\
Now we proceed to define a new Riemannian structure on the manifold so that the
action $\alpha$ will be inner product preserving. For that we are going to need
the following \\
\blmma
\label{new_Rieman}
(i) For $m\in M$, Sp $\{ds(m): s\in\cla\}$ coincides with $T^{*}_{m}(M)$.\\
(ii) If $\{ s_{1},...,s_{n}\}$ and $\{s_{1}^{'},...,s_{n}^{'}\}$ are two sets of
functions in $\cla$ such that each of $\{ds_{i}(m):i=1,...,n\}$ and
$\{ds_{i}^{'}(m):i=1,...,n\}$ are bases for $T_{m}^{\ast}(M)$ and for
$v,w\in T_{m}^{\ast}(M)$ with
$v=\sum_{i}c_{i}ds_{i}(m)=\sum_{i}c_{i}^{'}ds_{i}^{'}(m)$ and
$w=\sum_{i}d_{i}ds_{i}(m)=\sum_{i}d_{i}^{'}ds_{i}^{'}(m)$, then 
$$\sum_{i,j}\bar{c_{i}}d_{j}<<ds_{i},ds_{j}>>^{'}(m)=\sum_{i,j}\bar{c_{i}^{'}}d_
{ j }^{'} <<ds_{i}^{'},ds_{j}^{'}>>^{'}(m),$$
where $<<,>>^{'}$ is the new $C^{\infty}(M)$ valued inner product introduced
earlier.
\elmma
{\it Proof}: \\
Choosing a coordinate neighbourhood $U$ around $m$ and a set of coordinates
$x_{1},...,x_{n}$ we have $ds(m)= \sum_{i=1}^{n}\frac{\partial s}{\partial
x_{i}}(m)dx_{i}(m)$.\\
Pick any $\eta\in T^{*}_{m}(M)$ i.e. we have $\eta=
\sum_{i=1}^{n}c_{i}dx_{i}(m)$ for some $c_{i}$'s in $\mathbb{R}$.\\
Choose any $f\in C^{\infty}(M)$ with $\frac{\partial f}{\partial x_{i}}(m)=
c_{i}$.\\
For $f\in C^{\infty}(M)$, by Fr\'echet density of $\cla$ we have  a sequence
$s_{n}\in \cla$ and an $n_{0}\in \mathbb{N}$ such that $$|\frac{\partial
s}{\partial x_{i}}(m)-\frac{\partial f}{\partial x_{i}}(m)|< \epsilon \ \forall
\ n\geq n_{0}.$$
So Sp $\{ds(m); s\in \cla\}$ is dense in $T^{*}_{m}(M)$. $T^{*}_{m}(M)$ being
finite dimensional Sp $\{ds(m): s\in\cla\}$ coincides with $T^{*}_{m}(M)$. Which
proves (i).  \\
\indent For proving (ii) first we prove the following fact:\\
{\it Let $m\in M$ and $\omega \in \Omega^{1}(\cla)$ such that $\omega= 0$ in a
neighbourhood $U$ of $m$. Then $<<\omega, \eta>>^{'}= 0$ for all $\eta\in
\Omega^{1}(\cla)$}\\
\indent For the proof of the above fact Let $V\subset U$ such that $V\subset
\bar{V}\subset U$.\\
Choose $f\in C^{\infty}(M)_{\mathbb{R}}$ such that $supp(f)\subset \bar{V}$,
$f\equiv 1$ on $V$ and $f\equiv 0$ outside $U$.\\
So we can write $\omega=(1-f)\omega$. Then
\begin{eqnarray*}
&&<<\omega,\eta>>^{'}(m)\\
&=& <<(1-f)\omega, \eta>>^{'}(m)\\
&=& <<\omega,\eta>>^{'}(m)(1-f)(m) \ (by \ (\ref{1111}))\\
&=& 0.
\end{eqnarray*} 
\indent Applying the above fact we can show:\\
{\it Let $m\in M$ and $\omega=\omega^{'}, \ \eta=\eta^{'}$ in a neighbourhood
$U$ of $m$. Then $<<\omega,\eta>>^{'}= <<\omega^{'},\eta^{'}>>^{'}$, $\forall
\omega,\omega^{'},\eta, \eta^{'}\in \Omega^{1}(\cla)$.}\\
\indent For the proof it is enough to observe that $<<\omega,\eta>>^{'}(m)-
<<\omega^{'},\eta^{'}>>^{'}(m)=<<\omega-\omega^{'},\eta>>^{'}(m)+
<<\omega^{'},\eta-\eta^{'}>>(m)$.\\
\indent As $\{ds_{1}(m),..., ds_{n}(m)\}$ and $\{ds^{'}_{1}(m),...,
ds^{'}_{n}(m)\}$ are two bases
for $T^{*}_{m}(M)$. Then they are actually bases for $T^{*}_{x}(M)$ for $x$ in a
neighbourhood $U$ of $m$. So there are $\{f_{ij}:i,j= 1(1)n\}$ in
$C^{\infty}(M)$ such that $$ds_{i}= \sum _{j=1}^{n}f_{ij}ds^{'}_{j}$$ on $U$ for
all $i= 1,\ldots,n$. Hence by the previous discussion
\begin{eqnarray}
<<ds_{i},ds_{j}>>^{'}(m)=<<\sum_{k}f_{ik}ds^{'}_{k},
\sum_{l}f_{jl}ds^{'}_{l}>>^{'}(m)
\end{eqnarray}
Let $v=\sum_{i=1}^{n}c_{i}ds_{i}(m)=\sum_{i=1}^{n}c^{'}_{i}ds^{'}_{i}(m)$ and
$w=\sum_{i=1}^{n}d_{i}ds_{i}(m)=\sum_{i=1}^{n}d^{'}_{i}ds^{'}_{i}(m)$. So by
definition
\begin{eqnarray*}
<v,w>^{'}&=& \sum_{ij}\bar{c_{i}}d_{j}<<ds_{i},ds_{j}>>^{'}(m)\\
&=&
\sum_{ijkl}\bar{c_{i}}d_{j}\bar{f_{ik}}(m)f_{jl}(m)<<ds^{'}_{k},ds^{'}_{l}>>^{'}
(m) \ ( \ by \ (\ref{1111})\\
&=& \sum_{kl}\bar{c_{k}}^{'}d_{l}^{'}<<ds_{k}^{'},ds_{l}^{'}>>^{'}(m)
\end{eqnarray*}
\qed\\
\indent {\it Proof of Theorem \ref{average}}:\\
Now we can define a new inner product on the manifold $M$. For
that let $v,w\in T^{\ast}_{m}(M)$ by (i) of Lemma \ref{new_Rieman} we choose
$s_{1},...,s_{n}\in\cla$ such that $ds_{1}(m),...,ds_{n}(m)$ is a basis for
$T^{\ast}_{m}(M)$. Let $\{c_{i},d_{i}:i=1,...,n\}$ be such that $v=\sum_{i}
c_{i}ds_{i}(m)$ and $w=\sum_{i}d_{i}ds_{i}$. Then we define
$$<v,w>^{'}:= \sum_{i,j}\bar{c_{i}}d_{j}<<ds_{i},ds_{j}>>^{'}(m).$$
It is evident that this is a semi definite inner product. We have to show that
this is a positive definite inner product. To that end let $<v,v>^{'}=0$ i.e. 
$$\sum_{i,j}\bar{c_{i}}c_{j}<<ds_{i},ds_{j}>>^{'}(x)=0,$$
where $v=\sum_{i}c_{i}ds_{i}(x)\in T_{x}^{\ast}(M)$. Since the Haar state $h$ is
faithful on $\clq_{0}$ and by assumption $<<d\alpha(ds_{i}),d\alpha(ds_{j})>>\in
\cla\ot\clq_{0}$, we can deduce that 
$$\sum_{i,j}\bar{c_{i}}c_{j}((id\ot m)(id\ot\kappa\ot id)(\alpha\ot
id)<<d\alpha(ds_{i}),d\alpha(ds_{j})>>)(x)=0.$$
Since $\epsilon \circ \kappa=\epsilon$ on $\clq_{0}$, applying
$\epsilon $ to the above equation, we get
$$ \sum_{i,j}\bar{c_{i}}c_{j}((id\ot m)(id\ot\epsilon\ot\epsilon)(\alpha\ot
id)<<d\alpha(ds_{i}),d\alpha (ds_{j})>>)(x)=0.$$ Using the fact that $\epsilon$ is $\ast$-homomorphism we get
$$
\sum_{i,j}\bar{c_{i}}c_{j}<({\rm id} \ot \epsilon)(d\alpha(ds_{i}))(x),(\rm id \ot \epsilon)(d\alpha(ds_{j}))(x)>=0.$$ 
It is easy to see that
$({\rm id} \ot \epsilon)(d\alpha(ds_{i}))=ds_{i}$ for all $i$. Hence we conclude that
$$<\sum_{i}c_{i}ds_{i}(x),\sum_{i}c_{i}ds_{i}(x)>=0,$$
i.e $<v,v>=0$ and hence $v=0$ (as $<\cdot, \cdot>$ is strictly positive definite, being an inner product on $T_x^*M$)
so that $<\cdot, \cdot>^\prime$ is indeed strictly positive definite, i.e. inner product. 
\indent We have already noted ( (ii) of Lemma \ref{new_Rieman})   that our definition is independent
of choice of $s_{i}$'s, and also  that with respect to this new
Riemannian structure on the manifold,
$\alpha$ is inner product preserving. This completes the proof of 
the Theorem \ref{average} on $\Omega^{1}(\cla)$ and hence on
$\Omega^{1}(C^{\infty}(M))$.
\qed\\

\end{document}